\documentclass[11pt,a4paper]{article}
\usepackage{cite}
\usepackage{amsmath, amsthm,mathtools,url ,amssymb,upgreek,slashed}
\usepackage{ifpdf}
\ifpdf
  \usepackage[pdftex]{graphicx}
  \usepackage{epstopdf}
\else
  \usepackage[dvips]{graphicx}
\fi
\parskip=\baselineskip
\def\Bbb{\mathbb}
\def\Tr{{\rm Tr}}
\def\16{{\bf 16}}
\def\1{{\bf 1}}
\def\Sdiff{{\mathrm{Sdiff}}}
\def\2{{\bf 2}}
\def\3{{\bf 3}}
\def\4{{\bf 4}}
\def\vol{{\mathrm{vol}}}
\def\sg{{\mathrm g}}
\def\i{{\mathrm i}}

\def\O{{\mathcal O}}

\def\be{\begin{equation}}
\def\ee{\end{equation}}

\def\R{{\Bbb{R}}}\def\Z{{\Bbb{Z}}}

\def\hat{\widehat}
\def\Pf{{\mathrm{Pf}}}
\def\D{{\mathcal D}}
\font\teneurm=eurm10 \font\seveneurm=eurm7 \font\fiveeurm=eurm5
\newfam\eurmfam
\textfont\eurmfam=\teneurm \scriptfont\eurmfam=\seveneurm
\scriptscriptfont\eurmfam=\fiveeurm

\font\teneusm=eusm10 \font\seveneusm=eusm7 \font\fiveeusm=eusm5
\newfam\eusmfam
\textfont\eusmfam=\teneusm \scriptfont\eusmfam=\seveneusm
\scriptscriptfont\eusmfam=\fiveeusm

\font\tencmmib=cmmib10 \skewchar\tencmmib='177
\font\sevencmmib=cmmib7 \skewchar\sevencmmib='177
\font\fivecmmib=cmmib5 \skewchar\fivecmmib='177
\newfam\cmmibfam
\textfont\cmmibfam=\tencmmib \scriptfont\cmmibfam=\sevencmmib
\scriptscriptfont\cmmibfam=\fivecmmib

\numberwithin{equation}{section}

\def\d{\mathrm d}

\def\SL{{\mathrm{SL}}}
\def\PSL{{\mathrm{PSL}}}
\def\cU{{\mathrm U}}
\def\Z{{\Bbb Z}}

\def\M{{\mathcal M}}
\def\fH{{\mathfrak H}}
\def\fM{{\mathfrak M}}

\def\M{{\mathcal M}}

\def\d{{\mathrm d}}

\def\J{{\sf H}}
\def\H{{\mathcal H}}
\def\D{{\mathcal D}}
\def\Bbb{\mathbb}

\def\diag{{\mathrm{diag}}}
\def\Tr{{\mathrm{Tr}}}
\def\R{{\Bbb R}}

\def\U{{\mathrm U}}
\def\Z{{\Bbb Z}}
\def\hh{{\sf K}}
\def\HH{{\sf H}}

\def\frak{\mathfrak}
\def\MM{{\frak M}}
\def\i{{\mathrm i}}
\def\PSL{{\mathrm{PSL}}}
\def\OSp{{\mathrm{OSp}}}
\def\Pf{{\mathrm{Pf}}}
\def\diff{{\mathrm{diff}}}
\def\SL{{\mathrm{SL}}}
\def\Ber{{\mathrm{Ber}}}
\def\h{\widehat}
\def\osp{{\mathrm{osp}}}

\def\X{{\mathcal X}}
\def\sv{{\sf V}}
\def\Y{{\mathcal Y}}

\begin{document}
\begin{titlepage}
\begin{flushright}

\end{flushright}

\vskip 1.5in
\begin{center}
{\bf\Large{Volumes And Random Matrices}}
\vskip
0.5cm { Edward Witten} \vskip 0.05in {\small{ \textit{Institute for Advanced Study}\vskip -.4cm
{\textit{Einstein Drive, Princeton, NJ 08540 USA}}}
}
\end{center}
\vskip 0.5in
\baselineskip 16pt
\begin{abstract}   This article is an introduction to newly discovered relations between volumes of moduli spaces of Riemann surfaces or super Riemann surfaces,
simple models of gravity or supergravity in two dimensions, and random matrix ensembles.     (The article is based
 on a lecture at the conference on the Mathematics of Gauge Theory and
String Theory, University of Auckland, January 2020.  It has been submitted to a special issue of the Quarterly Journal of Mathematics in memory
of Michael Atiyah.)

 \end{abstract}
\date{March, 2020}
\end{titlepage}
\def\Hom{\mathrm{Hom}}

\def\U{{\mathrm U}}
\def\calU{{\mathcal U}}

\tableofcontents
\section{Introduction}\label{intro}

In this article, I will sketch some recent developments involving the volume of the moduli space $\M_\sg$ of Riemann surfaces of genus $\sg$, and also the volume of the corresponding
moduli space $\fM_\sg$ of super Riemann surfaces.   The discussion also encompasses Riemann surfaces with punctures and/or boundaries.   The main goal
is to explain how these volumes are related to random matrix ensembles.   This is an old story with a very contemporary twist.  

What is meant by the volume of $\M_\sg$?   One answer is that $\M_\sg$ for $\sg>1$ has a natural Weil-Petersson symplectic form.   A Riemann surface $\Sigma$ of genus $\sg>1$ can be
regarded as the quotient of the upper half plane $\fH\cong \PSL(2,\R)/\cU(1)$ by a discrete group.   Accordingly $\Sigma$ carries a natural Riemannian metric $g$ pulled
back from $\H$; this metric has 
constant scalar curvature $R=-2$ and is called a hyperbolic metric.   Similarly,  $\Sigma$ is endowed with a natural flat $\PSL(2,\R)$ bundle\footnote{\label{notee} The group $\PSL(2,\R)$ is contractible onto $\cU(1)$, so a $\PSL(2,\R)$ bundle over a surface $\Sigma$ has an integer invariant, the first Chern class.
The flat bundle related to a hyperbolic metric has first Chern class $2-2\sg$.   A similar remark applies for the supergroup $\OSp(1|2)$ introduced shortly.
The maximal bosonic subgroup of $\OSp(1|2)$ is the spin double cover $\SL(2,\R)$ of $\PSL(2,\R)$, so in that case the relevant value of the first Chern class is $1-\sg$.} with 
connection $A$ pulled back from $\H$. A point in $\M_\sg$ determines a Riemann surface together with such a flat connection, and the symplectic form $\omega$ of $\M_\sg$ can be defined by 
\be\label{sympform}\omega =\frac{1}{4\pi}\int_\Sigma \Tr \,\delta A\wedge \delta A ,    \ee 
in close analogy with the definition 
used by Atiyah and Bott  \cite{AB} for a symplectic form on the moduli space of flat bundles with compact structure group $G$.
The volume of $\M_\sg$ is then
\be\label{vol}V_\sg=\int_{\M_\sg}\Pf(\omega) = \int_{\M_\sg}e^\omega,\ee
where $\Pf$ is the Pfaffian.

This approach generalizes perfectly well for super Riemann surfaces.   One replaces $\SL(2,\R)$, which is the group of linear transformations of $\R^2$ that preserve
the  symplectic form $\d u\,\d v$, with $\OSp(1|2)$, which is the  supergroup  of linear transformations of $\R^{2|1}$ that preserve the symplectic
form $\d u\d v-\d\theta^2$. $\OSp(1|2)$ is a Lie supergroup of dimension $3|2$.   
Its Lie algebra carries a nondegenerate bilinear form that
I will denote as $\Tr$.  
 A point in $\MM_\sg$, the moduli space of super 
 Riemann surfaces of genus $\sg$, determines an ordinary Riemann surface $\Sigma$ together with a flat $\OSp(1|2)$ connection $A$.
 The symplectic form of $\MM_\sg$ can be defined by the same formula as before, only for flat $\OSp(1|2)$ connections rather than flat $\PSL(2,\R)$ connections:
 \be\label{superform}\hat\omega =\frac{1}{4\pi}\int_\Sigma \Tr \,\delta A\wedge \delta A ,    \ee 
The volume $\h V_\sg$ of $\MM_\sg$ can be defined 
\be\label{wory}\h V_\sg=\int_{\MM_{\sg}}\sqrt{\Ber\,\h\omega},\ee
where $\Ber$ is the Berezinian, the superanalog of the determinant.  This is analogous
to the first of the two formulas in eqn. (\ref{vol}).  There is no good superanalog of the integral of $e^\omega$, but the formula $\Pf(\omega)$ for a measure
does have a good superanalog, namely $\sqrt{\Ber(\h\omega)}$.  

I will now make a slight digression to explain how an ordinary Riemann surface $\Sigma$ with a flat $\OSp(1|2)$ connection  of the appropriate 
topological type (see footnote \ref{notee})
determines a super Riemann surface $\h\Sigma$.    
The superanalog of the upper half plane $\H$ is $\h \H=\OSp(1|2)/\U(1)$.    $\h \H$ is a smooth supermanifold of real dimension $2|2$; it
carries a complex structure in which it has complex dimension $1|1$.  
$\h \H$ also carries a 
canonical ``completely unintegrable distribution'' making it
a super Riemann surface.   
(There is no natural splitting of the Lie superalgebra $\osp(1|2)$ as the direct sum of even and odd parts, but the choice of a point in $\h \H$
determines such a splitting, and the odd part defines a subbundle of the tangent bundle to $\h \H$.  This is the unintegrable distribution.)
By taking the monodromies of the flat connection $A\to\Sigma$, we get a homomorphism $\rho:\pi_1(\Sigma)\to \OSp(1|2)$.
Let $\Gamma=\rho(\pi_1(\Sigma))$.   The quotient $\h\Sigma=\h\H/\Gamma$ is a smooth supermanifold of dimension $2|2$ that inherits from $\h\H$ the structure
of a super Riemann surface.  $\h\Sigma$ is the super Riemann surface associated to the pair $\Sigma,A$.   For our purposes here, however, it is convenient
to study the pair $\Sigma,A$ rather than the super Riemann surface $\h\Sigma$.  

 It is possible to describe the super volumes $\h V_\sg$ in purely bosonic terms, that is in terms of ordinary geometry.
 The ``reduced space'' of $\MM_\sg$ is the moduli space $\M'_\sg$ that parametrizes an ordinary Riemann surface $\Sigma$
with a spin structure, which we can think of as a square root $K^{1/2}$ of the canonical bundle $K\to \Sigma$.  $\M'_\sg$ is a finite cover of $\M_\sg$.
 The normal bundle to $\M'_\sg$ in $\MM_\sg$ is the vector bundle $U\to\M'_\sg$ whose fiber is $H^1(\Sigma,K^{-1/2})$. 
Viewing $U$ as a real vector bundle (of twice its complex dimension), we denote its Euler class as $\chi(U)$.  
The symplectic form $\h \omega$ of $\MM_\sg$ restricts along $\M'_\sg$ to the ordinary symplectic form $\omega$ of $\M'_\sg$ (which is a finite
cover of $\M_\sg$). 
 By general arguments about symplectic supermanifolds, one can show that
 \be\label{canform}\h V_\sg=\int_{\M'_\sg}\chi(U) e^{\omega}. \ee
 Thus what I will say about  supervolumes can be interpreted as a purely classical statement about $\M'_\sg$.   The quantities on the right hand side of eqn. (\ref{canform})
 were studied by Norbury \cite{Norbury} from a different point of view (related however to the same spectral curve that appears in our discussion in section \ref{super}).  
 
 The bosonic volumes $V_\sg$ can also be expressed in terms of the intersection theory of $\M_\sg$.   Indeed, the Weil-Petersson form is one of the tautological classes on
 $\M_\sg$ that were
 introduced by Mumford, Morita, and Miller.    Maryam  Mirzakhani
 described  \cite{M2} a sort of converse of this statement: from a knowledge of the volumes (for Riemann surfaces possibly with geodesic boundary, as discussed momentarily)
 one can deduce the intersection numbers of tautological classes.    These facts do not generalize directly to $\MM_\sg$, as there is not a natural intersection theory on a supermanifold.
 
  Volumes for surfaces with boundary are introduced as follows.    Let $\Sigma$
be a hyperbolic Riemann surface of genus $\sg$ with $n$ boundaries.   We require the boundaries to be geodesics of prescribed lengths $b_1,b_2,\cdots,b_n$.
  Let $\M_{\sg,\vec b}$ be the moduli space of such objects.   It has a symplectic form and volume that can be defined
  by precisely the same formulas (\ref{sympform}) and (\ref{vol}) as before.  
  Mirzakhani showed that $V_{\sg,\vec b}$ is a polynomial in $b_1,b_2,\cdots,b_n$, and that the canonical intersection numbers are the coefficients
of the top degree terms in this polynomial.  

One can similarly define a moduli space $\MM_{\sg,\vec b}$ of super Riemann surfaces with geodesic  boundaries of specified lengths $b_1,b_2,\cdots, b_n$, and its volume $\h V_{\sg,\vec b}$.

 The relation of volumes to intersection numbers gives one way to compute them \cite{Z,MZ}.
Mirzakhani, however, discovered a new direct way to compute the volumes \cite{M1}.    
  Regardless of how volumes are to be computed, the relation between volumes and intersection numbers shows
 that volumes are related to random matrix ensembles, since intersection numbers are known
 to be related to random matrix ensembles in multiple ways.
    Indeed, my conjecture \cite{WittenOld} on intersection theory on $\M_\sg$ was inspired directly
 by discoveries made in that period \cite{DS,GM,BK} relating two dimensional gravity to random matrix ensembles.  Moreover, Kontsevich's proof \cite{K} involved
a relation of the intersection numbers to a different random matrix ensemble.   However, at least in my opinion, the role of the random matrices in all of these
constructions was somewhat obscure.   The new developments that I will describe give a much clearer picture, since the random matrix will have a simple physical meaning.

In this article, I will sketch how Saad, Shenker, and Stanford   \cite{SSS}, following Eynard and Orantin \cite{EO}, reinterpreted Mirzakhani's results in terms of a random
matrix ensemble.  They  were motivated by considerations of quantum gravity, and I will try to give at least some idea about their motivation and
why their results are of physical interest.      Then 
I will explain how Stanford and I \cite{SW} generalized the story   to super Riemann surfaces and  quantum supergravity.  

\section{Universal Teichm\"{u}ller Space}\label{uniteich}  

Let $S^1$ be a circle and $\diff\,S^1$ its group of orientation-preserving diffeomorphisms.  The homogeneous spaces  $\diff S^1/\PSL(2,\R)$ or $\diff S^1/\cU(1)$ can be 
viewed as coadjoint orbits
of a central extension of
$\diff\,S^1$, so they carry natural symplectic structures.   Let $\X$ be one of these spaces, and denote its symplectic structure by $\omega$.  We will see that $\X$ with
its symplectic structure is a sort of infinite-dimensional analog of $\M_\sg$.   Actually,
 $\diff S^1/\PSL(2,\R)$ is sometimes
called universal Teichm\"{u}ller space, and the developments that will be reviewed here  perhaps give a perspective on the sense in which this name is justified.

It is believed that there is no reasonable definition of the ``volume'' $\int_\X e^\omega$.  However, we can do the following.  
Consider a subgroup $\cU(1)\cong S^1\subset \diff S^1$, consisting of rigid rotations of $S^1$.   In other words,
for some parametrization of $S^1$ by an angle $\theta$, $\cU(1)$ acts by $\theta\to\theta+{\mathrm{constant}}$.   
Then there is a moment map $\hh$ for this action of $\cU(1)$; in other words, if $\sv$ is the vector field on $\X$ that generates $\cU(1)$ and $\i_\sv$ is contraction with
$\sv$, then
\be\label{mmap} \d \hh=-\i_\sv \omega.\ee
Then introducing a real constant $\beta$, the integral
\be\label{wokko}Z(\beta) =\int_\X\exp(\hh/\beta+\omega)  \ee
is better-behaved.

In fact, this integral can be viewed as an infinite-dimensional example of a situation that was studied by Duistermaat and Heckman \cite{DH}, 
and then reinterpreted by Atiyah and Bott \cite{AB2}  in terms of
equivariant cohomology.
 Let $\Y$ be a symplectic manifold (compact, or with some suitable conditions at infinity) with symplectic form $\omega$ and action of $\cU(1)$.     Let $p_1,\dots,p_s$ be
the fixed points of the $\U(1)$ action.   For simplicity I assume that there are finitely many.    Let $\hh$ be the moment map for the $\U(1)$ action.
 The Duistermaat-Heckman/Atiyah-Bott (DH/AB) formula gives
$$\int_\Y  \exp(\hh/\beta+\omega)= \sum_i \frac{\exp(\hh(p_i)/\beta)}{\prod_\alpha( e_{i,\alpha}/2\pi \beta) },$$
 where the $e_{i,\alpha}$ are integers that represent the eigenvalues of the $\U(1)$ action on the tangent space to $\Y$ at $p_i$.
 
In the present example, there is only one fixed point in the $\U(1)$ action on $\diff S^1/\PSL(2,\R)$ or $\diff S^1/\U(1)$.    
The product over eigenvalues at this fixed point becomes formally $\prod_{n=2}^\infty  n/2\pi\beta$  in the example of $\diff S^1/\PSL(2,\R)$
(or  $\prod_{n=1}^\infty \,n/2\pi\beta$ in the other example). 
This infinite product is  treated with (for example) $\zeta$-function regularization. For $\diff S^1/\PSL(2,\R)$, the result is
\be\label{poly} Z(\beta) =\frac{C}{4 \pi^{1/2} \beta^{3/2}} \exp(\pi^2/\beta),\ee
where the constant $C$, which has been normalized for later convenience, depends on the regularization and so is considered inessential, but the rest is ``universal.''
(This  problem was first studied by  Kitaev \cite{Ki} and subsequently analyzed in \cite{MS}.   There are many derivations of the formula (\ref{poly}) in the physics literature.   The explanation that
I have sketched is from \cite{SW2}.)  There is a similar formula for the other example $\diff S^1/\U(1)$.   

 To use the DH/AB  formula, we did not need to know what is the moment map $\hh$ (only its value at the fixed point). 
  But in fact it is a function of interest.    To pick the $\U(1)$ subgroup of $\diff \,S^1$ that was used in this ``localization,'' we had to
pick an angular parameter $\theta$ on the circle; an element of $\diff\,S^1$ maps this to another parameter $u$, and $\hh$ is the integral of the Schwarzian
derivative $\{u,\theta\}$. 

By an inverse Laplace transform, we can write 
\be\label{zub} Z(\beta)=\int_0^\infty \d E \rho(E)\exp(-\beta E),\ee
with
\be\label{plub} \rho(E)= \frac{C}{4\pi^2} \sinh(2\pi\sqrt  E).  \ee

\section{Quantum Gravity In Two Dimensions}\label{qgtwo}

I will next explain why this result was considered problematical and how it has been interpreted \cite{SSS}.  This will require explaining some physics.

 General Relativity is difficult to understand as a quantum 
 theory.    Searching for understanding, physicists have looked for a simpler model in a lower dimension.    Two dimensions
 is a good place to look.      An
 obvious idea might be to start with the Einstein-Hilbert action in two dimensions, $I_{EH}=\int_\Sigma \d^2x \sqrt g R/2\pi$, with $R$ the
 Ricci scalar of a Riemannian metric $g$.   This does not work well, as in two-dimensions this action is a topological
 invariant, the Euler characteristic $\chi$, according to the Gauss-Bonnet theorem.    Instead it turns out to be better to add a scalar (real-valued) field
 $\phi$.     For many purposes, a simple and illuminating model of two-dimensional gravity is ``Jackiw-Teitelboim (JT) gravity,''
 with action
 \be\label{polfo} I_{JT}=-\frac{1}{2}\int_\Sigma\d^2x\sqrt g\, \phi(R+2). \ee
 Actually, even though the Einstein-Hilbert action is a topological invariant,
 it turns out that it is important to include a multiple of this term in the action.   The combined action is then
 \be\label{noffo}I=I_{JT} -S_0 I_{EH} =-\frac{1}{2}\int_\Sigma\d^2x\sqrt g\, \phi(R+2)- S_0 \int \d^2s\sqrt g \frac{R}{2\pi}. \ee
In conventional General Relativity, the coefficient of the Einstein-Hilbert term is $1/16\pi G_N$, where $G_N$ is Newton's constant.   
In the real world,  $G_N$ is extremely small in natural units set by the values of Planck's constant, the speed of light, and atomic masses.   
So $1/16\pi G_N$ is very large in natural units.  Similarly, here we are going to think of $S_0$ as being large.  The same goes for
the renormalized parameter $S$ introduced below.

Because $I_{EH}$ is a topological invariant, the Euler-Lagrange equations come entirely from $I_{JT}$.
The Euler-Lagrange equation for $\phi$ is simply $R+2=0$, so a classical solution is a hyperbolic Riemann surface.      The Feynman path integral for
  compact $\Sigma$ without boundary (or with geodesic boundary of prescribed length)   is very simple.    The path integral is
  \be\label{woodo} Z_\Sigma=\frac{1}{\vol}\int\D \phi\,\D g \,e^{-I},\ee
  where formally $\vol$ is the volume of the diffeomorphism group and the factor $1/\vol$ is a way to indicate that we have to divide by the diffeomorphism
  group (that is, we integrate over pairs $\phi,g$ up to diffeomorphism).   Since $I=I_{JT}-S_0 I_{EH}=I_{JT}-S_0\chi$, we have in more detail
  \be\label{wodo}Z_\Sigma= e^{S_0\chi}\cdot
   \frac{1}{\mathrm{vol}}\int \D\phi\,\D g \exp\left(\frac{1}{2}\int\d^2x\sqrt g \phi(R+2)\right).\ee
  The integral over $\phi$ and $g$ 
  is studied by integrating first over $\phi$ (after rotating the integration contour $\phi\to \i\phi$) and gives a delta function setting $R+2=0$.
Since we have to divide by the diffeomorphism group, the integral 
  ``localizes'' on the moduli space of two-manifolds with hyperbolic structure, modulo diffeomorphism.
  
  If $\Sigma$ is orientable and of genus $\sg$, the moduli space of two-manifolds with metric of constant curvature $R=-2$ is the usual moduli space $\M_\sg$
  of Riemann surfaces of genus $\sg$, and one can show that the integral over $\M_\sg$ gives its usual volume, up to an elementary factor:
  \be\label{usvol} \frac{1}{\mathrm{vol}}\int \D\phi\,\D g \exp\left(\frac{1}{2}\int\d^2x\sqrt g \phi(R+2)\right) =  C^{\chi}\int_{\M_\sg} e^\omega.    \ee
$C$ is a constant, independent of $\sg$, that depends on the regularization used
  in defining the integral.   There is no natural choice of  regularization, so there is no preferred value of the constant $C$.
 However, in eqn. (\ref{woodo}), we see that $Z_\Sigma$ has an additional factor $e^{S_0\chi}$.  We simply eliminate the arbitrary constant $C$ by setting $e^S=e^{S_0}C$.
 So finally 
 \be\label{finalmy} Z_\Sigma = e^{S\chi} \int_{\M_\sg} e^\omega. \ee
 What has just been explained is a typical, though elementary, example of renormalization theory.
 The theory does not depend on the choice of regularization, as long as it is expressed in terms of the ``renormalized'' parameter $S$ rather than the ``bare''
 parameter $S_0$.     If $\Sigma$ is unorientable, the moduli space of hyperbolic structures is not a symplectic manifold,
 and the formula analogous to eqn. (\ref{finalmy})  involves the Reidemeister or Ray-Singer torsion \cite{SW}.

  So far, I have assumed that $\Sigma$ is a compact surface without boundary.  The case that $\Sigma$ has geodesic boundaries of specified length is similar.
       What really led to progress
  in the last few years was applying JT gravity to, roughly speaking, the whole upper half-plane $\H$.  
  The Euler-Lagrange  field equations of JT gravity tell us that $R+2=0$, and also give a certain equation for the real-valued field $\phi$.
  Since the natural metric on $\H$ has $R+2=0$, $\H$,
  endowed with a suitable $\phi$ field, can be regarded as a classical solution of JT gravity.   The Euler-Lagrange equations for $\phi$ have
  a natural interpretation if one views $\H$ as a Kahler manifold, and thus in particular as a Riemannian manifold that also has a symplectic
  structure.   These equations say that $\phi$ is the moment map for a vector field $V$ on $\H$ that generates an automorphism of $\H$
  as a Kahler manifold: thus $V$ is a Hamiltonian vector field, with Hamiltonian function $\phi$, for $\H$ viewed as a symplectic
  manifold, and also is a Killing vector field
  for $\H$ viewed as a Riemannian manifold.   In short, $V$ generates a one-parameter subgroup of the automorphism group
   $\PSL(2,\R)$ of $\H$.   The case that one wants (because it makes possible a limiting procedure that is described shortly)
   is that this is a compact subgroup $\U(1)\subset \PSL(2,\R)$.   Any two
   such subgroups are conjugate so it does not matter which one we choose.
   Concretely, if one describes $\H$ as the sheet $x>0$ in a hyperboloid $x^2-y^2-z^2=1$, then we can choose the $\U(1)$
   that rotates $y$ and $z$.   The moment map for a generator of this $\U(1)$ is simply the function $x$, and we take
  $\phi$ to be a positive multiple of this:     \be\label{wifgo} \phi=cx, ~~~~~c>0.\ee
  In this description, the conformal boundary of $\H$ is at $x\to\infty$, and we see that $\phi$ blows up everywhere near this
  conformal  boundary.

 But it turns out that literally
studying JT gravity on all of $\H$ is not the right thing to do.     This would be rather like trying to calculate the naive integral $\int_{\X} e^\omega$, with $\X=\diff S^1/\PSL(2,\R)$,
rather than the improved
  version (\ref{wokko}).   A better thing to do is to study JT gravity not on all of $H$ but on a very large region $U\subset \H$ \cite{MaldaEtal}.

  \begin{figure}
 \begin{center}
   \includegraphics[width=3.5in]{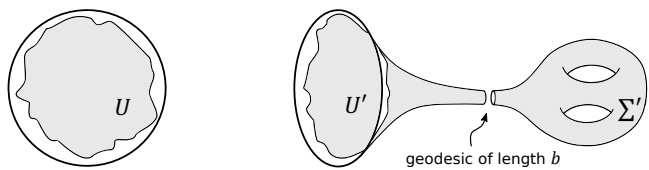}
 \end{center}
\caption{\small   On the left, the upper half plane $\H$ is represented as a disc, and the shaded portion of the disc is a large region $U\subset \H$.    On the right, a Riemann surface
$\Sigma$ is built by gluing together two pieces $U'$ and $\Sigma'$ (shaded) along an embedded geodesic of circumference $b$.   $\Sigma'$
 is a surface of $\sg>0$ with one geodesic
boundary.    $U'$   is an annulus with a geodesic ``inner''  boundary and a large ``outer'' boundary.    Near its outer boundary, $U'$ resembles locally a large region in $\H$.
\label{dt}}
\end{figure}
  
Such a large region is sketched on the left of fig. \ref{dt}.   (The right of the figure will be discussed later.)   
We are going to study JT gravity on a two-manifold $U$ that topologically
is a disc.  
To get a good variational problem for JT gravity on a manifold with boundary, one needs to include in the action a boundary contribution,
somewhat like the boundary correction to the Einstein-Hilbert action of 
 General Relativity
  \cite{York,GH}.  
The JT action, with the boundary term, is
\be\label{blab} I_{JT}=-\frac{1}{2}\int_\Sigma \d^2x\sqrt g \phi(R+2)-\frac{1}{2}\int_{\partial\Sigma}\d x\sqrt h \phi( K-1) , \ee
    where $K$ is the extrinsic curvature of the boundary $\partial \Sigma$, and $h$ is the induced metric of the boundary.   
    The Einstein-Hilbert action $I_{EH}$ similarly needs a boundary term, so that it still equals the Euler characteristic $\chi$. 
On $\partial U$, we impose a version of\footnote{Note in particular that we do {\it not} constrain $\partial U$ to be a geodesic.
Classically, it would be impossible to impose such a constraint, since a disc does not admit a hyperbolic metric with geodesic boundary.
But more to the point, we want a boundary condition that leads $U$ to be, in some sense, a good approximation to $\H$.}
 Dirichlet boundary conditions: we specify the induced metric $h$ of the boundary $\partial U$, along with the boundary value of\footnote{With this boundary
 condition, eqn. (\ref{blab}) leads to a sensible classical variational problem (and quantum path integral) regardless of whether we use $K$ or $K-1$
 in the boundary term.     We use $K-1$ because $K-1$  vanishes in the limit of a large disc in $\H$.
 This is important for the existence of the limit discussed momentarily in which $U$ becomes large.} $\phi$.
 Moreover, we specify them so that the circumference of $\partial U$ is very large, as is $\phi|_{\partial U}$.

 In a limit as the circumference of $\partial U$ goes to infinity and $\phi|_{\partial U}$ is taken to be a large constant, 
 the Feynman integral $Z_U(\beta)$ of JT gravity on $U$ can be evaluated in terms of an integral that we already studied in section \ref{uniteich}:
 $Z_U(\beta)=e^{S_0}\int_{\X}\exp(\hh/\beta+\omega)$, where $\X=\diff\,S^1/\PSL(2,\R)$.    The prefactor $e^{S_0}$ comes from the Einstein-Hilbert part of the action; that is, it is  $e^{S_0\chi}$, where the disc $U$ has $\chi=1$.    More subtle is to explain how the rest of the path integral gives 
 the integral $\int_{\X}\exp(\hh/\beta+\omega)$.     Note first that $S^1=\partial \H$ comes with an angular parameter that is uniquely determined up to the action
 of $\PSL(2,\R)$.   Moreover, since its induced metric is specified, $\partial U$ has a natural arclength parameter, up to an additive constant.   In the limit of interest in which $U$ is
 very large, it makes sense to compare these two parameters, and the comparison defines an element of $\X=\diff\,S^1/\PSL(2,\R)$.
  So we can think of $\X$ as parametrizing a family of large regions $U\subset \H$, up to $\PSL(2,\R)$.   
 It turns out   that in the limit that $\phi|_{\partial U}$ is large, forcing $\partial U$ closer to $\partial \H$,
 $\int_{\partial\Sigma}  (K-1)$ becomes a multiple of the moment map $\hh$ for the $\U(1)$ action on $\X$.   Hence in this limit, apart from the factor $e^{S_0}$, the Feynman integral $Z_U$
    of JT gravity on the disc $U$ becomes our friend 
    $ \int_{\X}\exp(\hh/\beta+\omega). $   To achieve convergence to this limit, both the circumference of $\partial U$ and the constant value of $\phi|_{\partial U}$ are taken
    to infinity, keeping fixed their ratio, $\beta$.           From eqns. (\ref{zub}) and (\ref{plub}), and setting $e^{S_0}C=e^S$, we get 
    \be\label{hubb} Z_U(\beta)=\int_0^\infty \d E\, \rho(E)\exp(-\beta E),~~\rho(E)=  \frac{e^S}{4\pi^2}\sinh(2\pi\sqrt  E).  \ee
Here  $\beta$ is regarded as the renormalized circumference of $\partial U$, and $S$ is the renormalized coefficient of the Einstein-Hilbert term.

   This is a deeply problematic answer for the Feynman integral on the disc $U$.
     To understand why, one should be familiar with holographic duality between gravity in the bulk of spacetime and an ordinary quantum system on the boundary \cite{Malda}.
      If the bulk were 4-dimensional, as in ordinary physics,
       the boundary would be 3-dimensional and the ``ordinary quantum system'' on the boundary would be a quantum field theory -- not a very
      easy concept to understand mathematically.
       But here the bulk is 2-dimensional, so the boundary is  1-dimensional and matters are simpler.
An ordinary quantum system in 1 dimension is just described by giving a Hilbert space $\J$ and a Hamiltonian operator $H$ acting on $\J$.
 The basic recipe of holographic duality predicts that $Z_U(\beta)=\Tr_\J \exp(-\beta H)$. 
    
    In a moment, we will check that that prediction is false, but before doing so, I want to explain that this actually did not come as a complete surprise.
    Analogous calculations going back to the 1970's \cite{Hawking,GH}  have always given the same problem:
   Euclidean path integrals give results that lack the expected Hilbert space interpretation.   
     The calculations were traditionally done in models (like four-dimensional General Relativity) that were too complicated for a complete calculation,
    and there was always a possibility that a more complete calculation would make the issue go away.   The novelty is that
     holographic duality and a variety of other related developments have made it possible to ask the question in a model -- JT gravity -- 
    that is so simple that one can do a complete calculation, demonstrating the problem.

    To see that the prediction of the duality is false, we just note the following.
     If we do have a Hilbert space $\HH$ and a Hamiltonian $H$ acting on it such that the operator $e^{-\beta H}$ has a trace, then $H$ must have a discrete
    spectrum with eigenvalues $E_1,E_2,\cdots$ (which moreover must tend to infinity fast enough) and
     \be\label{knop} \Tr_\HH\,\exp(-\beta H)=\sum_i e^{-\beta E_i} =\int_0^\infty \d E \,\sum_i\delta(E-E_i) \, e^{-\beta E}. \ee
      But the integral over $\diff\,S^1/\PSL(2,\R)$ gave
     \be\label{nop} Z_U(\beta)= \int_0^\infty \d E  \,\frac{e^S}{4\pi^2}\sinh(2\pi\sqrt  E) e^{-\beta E}. \ee
      The function $\frac{e^S}{4\pi^2}\sinh(2\pi\sqrt  E)$ is not a sum of delta functions, so the prediction of the duality is false.

 However, the interpretation via JT gravity gives us a key insight that we did not have when we were just abstractly integrating over  $\X=\diff S^1/\PSL(2,\R)$.   
 In equation (\ref{plub}), $C$ was an arbitrary constant, but now $C$ has been replaced with $e^S$ where $S$ is should be large to get a model more similar
 to the real world.   This was explained following eqn. (\ref{noffo}).   To amplify on the point a bit, $S$ is analogous to black hole entropy in four dimensions; it is
very large in any situation in which the usual puzzles of quantum black hole physics arise.
   If $S$ is just moderately large, say of order 100, then $e^S$ is huge.    So we  think of $e^S$ as a huge number.
 With this in mind, 
the function  $\frac{e^S}{4\pi^2}\sinh(2\pi\sqrt E)$ actually  can be well-approximated for many purposes as $\sum_i \delta(E-E_i)$, for suitable $E_i$.    One must look very closely to see the difference.     One could not reasonably 
 approximate $2 \sinh(2\pi\sqrt E)$ (for example) by a such a sum, but $e^{100}\sinh(2\pi\sqrt E)$   is another matter.  
   
    \section{A Random Matrix Ensemble}\label{rme} 
    
    The novel idea of Saad, Shenker, and Stanford \cite{SSS} was to interpret the function $\frac{e^S}{4\pi^2} \sinh 2\pi\sqrt{E}$   not as the density of energy levels of a particular
    Hamiltonian but as the average level density of an ensemble of Hamiltonians -- a random matrix.  In terms of the physics involved, this interpretation is rather provocative,
    though it may be a challenge to convey this fully.
    
    What motivated this interpretation?   One clue came from the work of Kitaev \cite{Ki}, who discovered a simple model of holographic duality based on a random ensemble
    that was more complicated than the one of Saad et. al.   Unfortunately, to explain this would take us rather far afield.  Another clue came from the prior history
    of relations between random matrices and two-dimensional gravity \cite{BK,DS,GM,K}.
    A final clue, more directly related to the topic of
    the present article, involved the volumes of the moduli spaces of Riemann surfaces.   As noted in the introduction, Mirzakhani \cite{M1} had found a powerful new
    way to compute the volumes of these moduli spaces. And Eynard and Orantin \cite{EO} had interpreted her formulas in terms of their notion of topological recursion \cite{EO2},
which is related to a random matrix ensemble.    The key technical observation of Saad et. al.  was that the eigenvalue density of a random matrix
ensemble related to the spectral curve of Eynard and Orantin is precisely the function $\frac{e^S}{4\pi^2}\sinh 2\pi\sqrt{E}$ that arises in JT gravity and integration over
$\diff\,S^1/\PSL(2,\R)$.

The  sort of random matrix ensemble used in \cite{SSS}  is the following.   $M$ will be an $N\times N$ hermitian matrix for some $N$.  We are really interested in $N$
very large and ultimately in a limit with $N\to\infty$.   Picking some suitable real-valued function $T(M)$, we consider 
the integral
\be\label{wollo} Z(T;N)=\frac{1}{\vol(\U(N))}\int \d M \exp(-N\Tr\,T(M)). \ee
(The space of $N\times N$ Hermitian matrices is a copy of $\R^{N^2}$, and the measure $\d M$ is the standard Euclidean measure on $\R^{N^2}$.) 
If the function $T(M)$ is quadratic, this is a Gaussian random matrix ensemble, as studied  originally in the 1960's by Wigner, Dyson, Mehta, and many others. 
We are interested in the case that $T(M)$ is not quadratic.   In this case, $Z(T;N)$ or more precisely its logarithm has an asymptotic expansion for large $N$:
\be\label{ollo}\log Z(T;N)=N^2F_0(T)+ F_1(T)+\frac{1}{N^2}F_2(T)+\cdots =\sum_{\sg=0}^\infty N^{2-2\sg} F_\sg(T). \ee
  This expansion can be constructed by standard Feynman diagram methods.   As shown by 
  't Hooft \cite{thooft},  the Feynman diagrams for a matrix model are in a natural way ribbon
  graphs, each of which can be naturally drawn on a certain Riemann surface.    $F_\sg(T)$ is
  the contribution of connected ribbon graphs that can be drawn on a surface of genus $\sg$.   
  
  However, instead of making a Feynman diagram expansion, we can try to just evaluate the integral \cite{bpiz}.   We do this by first diagonalizing $M$,
  writing $M=U \Lambda U^{-1}$, with $U\in \U(N)$, and $\Lambda$ a diagonal matrix $\Lambda=\diag(\lambda_1,\lambda_2,\cdots,\lambda_N)$, $\lambda_1\leq \lambda_2\leq\cdots
  \leq \lambda_N$.
  The measure becomes $\d M =\d U \prod_i\d\lambda_i \prod_{j<k}(\lambda_j-\lambda_k)^2$.   Here $\d U$ is Haar measure on the group $\U(N)$, and the integral over
  $U$ just cancels the factor $1/\vol(\U(N))$ in the definition of $Z(T;N)$.   So we reduce to
  \be\label{pollo}Z(T;N)=\int\d\lambda_1\cdots\d\lambda_N \prod_{i<j}(\lambda_i-\lambda_j)^2\exp\left(-N\sum_k T(\lambda_k)\right). \ee  
   
  \begin{figure}
 \begin{center}
   \includegraphics[width=2.5in]{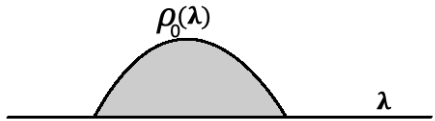}
 \end{center}
\caption{\small   For a nice class of functions $T(\lambda)$, the support of the  eigenvalue density $\rho_0(\lambda)$     is an interval $[a,b]$.
\label{Bump}}
\end{figure}

  For large $N$, the integrand $\prod_{i<j}(\lambda_i-\lambda_j)^2\exp(-N\sum_k T(\lambda_k))$ has a sharp maximum as a function of the $\lambda_i$.  To describe this maximum, assuming $N$ to be very large,  we think in terms of a continuous distribution of eigenvalues,  $\rho(\lambda)=N\rho_0(\lambda)$ for  some function $\rho_0(\lambda)$ (normalized to $\int\d\lambda \rho_0(\lambda)=1$).  In terms of such a distribution,
   the integrand in (\ref{pollo}) becomes
 \be\label{kolk}  \exp\left(N^2\left(-\int\d\lambda \rho_0(\lambda)T(\lambda) +\int\d\lambda\d\lambda'\rho_0(\lambda)\rho_0(\lambda') \log|\lambda-\lambda'|\right)\right). \ee
 For a nice class of functions $T(\lambda)$, the exponent has a unique maximum at some function $\rho_0(\lambda)$.     The most famous case is that $T$ is quadratic; then (with
 suitable normalization) $\rho_0(\lambda)=\frac{2}{\pi}\sqrt{1-\lambda^2}$ in the interval $[-1,1]$ and $\rho_0(\lambda)=0$ outside the interval.   This example illustrates a general behavior; for a nice
 class of functions $T$, the support of the function $\rho_0(\lambda)$  is an interval $[a,b]$; moreover, $\rho_0(\lambda)$ vanishes near the endpoints 
 of the interval as the square root of the distance to the endpoint, and
 can be continued to an algebraic function of $\lambda$ with only those singularities.\footnote{Because of the square root singularities at $\lambda=a,b $,
 the analytically continued function $\rho_0(\lambda)$
is imaginary on the real axis outside the interval $[a,b]$, so in terms of
 the continued function, the eigenvalue density is  $\mathrm{Re}\,\rho_0(\lambda)$.  For example, note that the function $\frac{2}{\pi}\sqrt{1-\lambda^2}$ which arises in the Gaussian
 case is imaginary on the real axis outside the interval $[-1,1]$.}

 If $F_0(T)$ is the value of the exponent at its maximum, then the leading approximation to the logarithm of the integral is
 \be\label{wallaby}\log Z(N;T)=N^2 F_0(T)+\O(1).  \ee   How can we calculate the further terms in this expansion?
 Something quite remarkable happens, though a detailed explanation would go beyond the scope of this article.  (There is a brief introduction in section 4 of 
 \cite{SW}, for example.)  
 One should define the ``spectral curve'' in the $y-\lambda$ plane:\footnote{The factor of $\pi^2$ here is conventional. Its purpose
 is to avoid a factor of $\pi$ in the relationship between $y$ and the expectation value of the matrix resolvent $\Tr\,\frac{1}{\lambda-M}$.}
 \be\label{wijo} y^2=-\pi^2\rho_0^2(\lambda). \ee
This curve is a double cover of the $\lambda$ plane, branched at the endpoints of the eigenvalue distribution, where the function $\rho_0^2(\lambda)$ has
simple zeroes.
 Once one knows this spectral curve, one can forget about doing integrals and one can forget the original function $T(M)$.   The whole
 expansion
 \be\label{nollo}\log Z(T;N)=N^2F_0(T)+ F_1(T)+\frac{1}{N^2}F_2(T)+\cdots =\sum_{\sg=0}^\infty N^{2-2\sg} F_\sg(T) \ee
(and a wide variety of other things about this ensemble) can be worked out just from a knowledge of the spectral curve.
A very useful version of this process is the ``topological recursion'' of Eynard and Orantin \cite{EO2}.   

Now let us go back to the volumes of moduli spaces.   As I have explained, Saad, Shenker, and Stanford interpreted the function $\frac{e^S}{4\pi^2} \sinh2\pi\sqrt{E}$ 
as the average density of eigenvalues for a random matrix drawn from a non-Gaussian hermitian matrix ensemble.   In principle, the procedure is to start with a function
$T$ and then compute the density of eigenvalues $N\rho_{0,T}(E)$, where we now denote the eigenvalue as $E$, and we
 also make explicit the dependence of 
$\rho_0$ on $T$.    Then we take $N\to\infty$ while adjusting $T$
so that $N\rho_{0,T}(E)$ converges to the desired $\frac{e^S}{4\pi^2} \sinh2\pi\sqrt{E}$.    This process is called double-scaling.  

But we can skip all that work, because everything we want to compute depends only on the spectral curve,  which we know is going to be
\be\label{tollo}y^2=-\frac{1}{16\pi^2}\sinh^2(2\pi\sqrt{E}). \ee
This curve is a double cover of the complex $E$-plane, with branch points at $E=0$ and $E=\infty$.   (In the process of double-scaling to get this curve,
the interval $[a,b]$ converges to $[0,\infty)$.) 
In short, all we have to do is to start with that spectral curve, and apply topological recursion to get the expansion
$\log Z(T)=\sum_{\sg=0}^\infty e^{(2-2\sg)S} F_\sg$, as well as other quantities of interest that are introduced momentarily.
Here, after double-scaling, the expansion parameter is $e^{-S}$ rather than $1/N$.

Now we can compute volumes, with a few steps, as follows.  First we compute the average of $\Tr\,\exp(-\beta H)$ in this
matrix ensemble, where now I refer to the random matrix $M$ as a Hamiltonian $H$.      This can be done explicitly, applying topological recursion
to the spectral curve (\ref{tollo}).   The result is an expansion of the ensemble average $\left\langle \Tr\,\exp(-\beta H)\right\rangle$ in powers of $e^{-2S}$. 

To interpret the result in terms of volumes, one proceeds as follows.   The Feynman diagram expansion of  $\left\langle
\Tr\,\exp(-\beta H)\right\rangle$ involves Feynman diagrams drawn on a Riemann surface $\Sigma$ with one boundary component, like those drawn in  fig. \ref{dt}.    
 This happens because when we make a Feynman diagram expansion, the trace $\Tr\exp(-\beta H)$ turns into
a boundary.  The picture on the left of fig. \ref{dt} actually corresponds to the leading contribution, the special case that $\Sigma$ is a disc.   In section \ref{qgtwo}, we
 analyzed the JT path integral for the case of a disc of regularized circumference $\beta$, 
and interpreted the answer as
$\int_0^\infty \d E\rho(E) \exp(-\beta E)$ with $\rho(E)=\frac{e^S}{4\pi^2}\sinh 2\pi\sqrt{E}$.   We have chosen a matrix ensemble that reproduces this
answer for the leading contribution.   We are now interested in the higher topologies shown on the right of fig. \ref{dt}.
They contribute the higher order terms in the expansion in powers of $e^{-S}$.   The contribution from Feynman diagrams drawn on a genus $\sg$ surface
$\Sigma$ is of relative order $\exp(-2\sg S)$.   

In the right hand panel of fig. \ref{dt}, the surface $\Sigma$ has been cut in two along a geodesic of length $b$.   To the right of the cut is a hyperbolic Riemann surface
$\Sigma'$ of genus $\sg$, with one geodesic boundary of length $b$. Let $\M_{\sg,b}$ be the moduli space of such hyperbolic
surfaces, and $V_{\sg,b}$ its volume.  To the left of the cut is a hyperbolic two-manifold $U'$ that is topologically an annulus.
$U'$ has one geodesic boundary  of length $b$ and another ``large'' boundary near which $U'$ 
 looks locally  like a large portion of $\H$.   In the same sense that the disc $U$ depicted on the left of fig. \ref{dt} represents $\diff S^1/\PSL(2,\R)$, $U'$ represents $\diff S^1/\U(1)$.  Apart from an overall multiplicative constant, the symplectic structure of  $\diff S^1/\U(1)$ depends
 on one parameter\footnote{As explained for example in section 3.1 of \cite{SW}, $\diff S^1/\U(1)$ can be viewed as a coadjoint orbit of the Virasoro group 
 (that is, the central extension of $\diff S^1$).   But actually, there is a one-parameter family of coadjoint orbits that are all isomorphic to $\diff S^1/\U(1)$, so the symplectic
 structure of $\diff S^1/\U(1)$ depends on a parameter $b$, in addition to the possibility of scaling it by an overall constant.} $b$, which corresponds to the length of the geodesic boundary of $U'$.
   Let us just write $\Theta(b;\beta)$ for the JT integral on this orbit ($\beta$ being again the regularized
circumference of the ``large'' boundary of $U'$).   
 The function $\Theta(b;\beta)$ can be computed from the DH/AB formula and is given by a formula similar to eqn. (\ref{zub}).
 Roughly because of locality properties of quantum field theory, or because the length of a geodesic on a hyperbolic two-manifold is the Hamiltonian
 that generates a corresponding Dehn twist,   the JT path integral $Z_\Sigma$ on $\Sigma$ 
 is obtained by integrating over $b$ the product of the JT path integral on $U'$ times the JT path integral on $\Sigma'$:
\be\label{zomo}Z_\Sigma=e^{S(1-2\sg)} \int_0^\infty b\d b \,\Theta(b;\beta) V_{\sg,b}.\ee
One has to include here a factor of $b$ that comes from integrating over a ``twist'' parameter in the gluing of $U'$ to $\Sigma'$.  We have used
the fact that the JT path integral on $\Sigma'$ is $e^{S\chi(\Sigma')}V_{\sg,b}=e^{S(1-2\sg)}V_{\sg,b}$.

On the other hand, the Feynman diagram expansion of the matrix model makes us expect that $Z_\Sigma$ will be the term of order
$e^{S(1-2\sg)}$ in the expansion of $\left\langle \Tr\,\exp(-\beta H)\right\rangle$.   This can be computed by applying topological recursion to the spectral curve (\ref{tollo}).
Using also the explicit result for $\Theta(b;\beta)$, one can make an inverse Laplace transform of
 the result for $Z_\Sigma$ to obtain an explicit formula for $V_{\sg,b}$.
 
 This procedure  can be generalized to give the volumes $V_{\sg,\vec b}$ of moduli spaces $\M_{\sg,\vec b}$ of hyperbolic surfaces of genus $\sg$
 with boundary geodesics of specified lengths $\vec b=(b_1,b_2,\cdots,b_n)$.   
The resulting formulas are correct, 
in view of the relation demonstrated by Eynard and Orantin \cite{EO} between topological recursion for the spectral curve $y^2=-\frac{1}{16\pi^2}\sinh^22\pi\sqrt{E}$ 
and the recursion relation used by Mirzakhani \cite{M1} to compute $V_{\sg,\vec b}$.

Matching with Mirzakhani's recursion relation was how Eynard and Orantin determined which spectral curve to use.   Another way is to
use the relation between volumes and intersection numbers and the general relation of intersection numbers to spectral curves.   (This approach is sketched in
\cite{DW}; see section 2.4 and eqn. (4.46).)  As I have explained, Saad, Shenker, and Stanford instead arrived at the same
spectral curve from the path integral of JT gravity on a disc.

The approach involving JT gravity is very interesting for physicists, but if one only cares about volumes of moduli spaces,
one might ask why it is important.  One answer is that this derivation sheds a new light on the relationship between $\diff\,S^1/\PSL(2,\R)$ and the moduli spaces
of Riemann surfaces.  Another answer is that this approach possibly gives a better understanding of why random matrix ensembles are related to volumes
and intersection numbers.  A third answer is given by my work with Stanford \cite{SW}.   We ran the entire
story for super Riemann surfaces.   Every step has a direct analog for that case.
 
 \section{Matrix Ensembles and Volumes of Super Moduli Spaces}\label{super}
 
 The superanalog of JT gravity is JT supergravity, which computes the volumes of moduli spaces of super Riemann surfaces,
 in general with geodesic boundaries of specified lengths.  As before, it is important to consider the special case of a super Riemann surface which is
 the super upper half plane $\h\H$, or more precisely a very large region $\h U\subset \h\H$, as in the left hand side of the familiar  fig. \ref{dt}.
 The boundary of $\h U$, or simply the boundary of $\h\H$, is a super circle that we may call $S^{1|1}$.   I also write $\Sdiff \,S^1$ for its group
 of orientation-preserving diffeomorphisms.   Let $\h\X=\Sdiff\,S^1/\OSp(1|2)$, which we might think of as universal super Teichm\"{u}ller space.   By similar arguments to those that were described in
 section \ref{qgtwo}, the super JT path integral on a large region $\h U\subset \h \H$ computes
 \be\label{woz} Z_{\h U}(\beta)=e^{S_0}\int_{\h\X}\exp(\h\hh/\beta+\h\omega), \ee
 where $\h\hh$ is the moment map for a subgroup $\U(1)\subset \Sdiff\,S^1$,  $\h\omega$ is the symplectic form of $\h\X$, and $S_0$ is the coefficient of the
 Einstein-Hilbert term in the action.
 This integral can be computed using DH/AB localization, with the result  \cite{SW} 
 \be\label{nz} Z_{\h U}(\beta)=\int_0^\infty \d E\, e^{-\beta E}\h\rho(E), ~~~\h\rho(E)= \frac{ e^S\,\sqrt 2 }{\pi}\frac{\cosh (2\pi\sqrt E)}{\sqrt E} .\ee
 As before, $S$ is the renormalized coefficient of the Einstein-Hilbert term.
 Again, this is not $\Tr_\HH\exp(-\beta H)$ for a Hamiltonian $H$ acting on a Hilbert space $\HH$, but now we know what to do: 
 we have to consider a random quantum ensemble.  
 
 We can rerun the previous story with a few changes.  The formula for $Z_{\h U}(\beta)$ tells us the spectral curve:
 \be\label{wz} y^2=-\frac{2}{E}\cosh^2(2\pi\sqrt E).  \ee
 However, the matrix ensemble is of a different type than we encountered before.
One way to see that it must be different is that the eigenvalue density $\h\rho(E)$ in eqn. (\ref{nz}) behaves as $E^{-1/2}$ near the endpoint of the energy spectrum at $E=0$,
in contrast to the typical $E^{+1/2}$ behavior of a matrix ensemble of the sort  that we studied in section \ref{rme}.

The ensemble must be different because as we are now studying super Riemann surfaces rather than ordinary ones,
   the dual quantum mechanical system is now supposed to be supersymmetric.
 So we need to do random supersymmetric quantum mechanics, not just random quantum mechanics.   
 
 In supersymmetric quantum mechanics, the Hilbert space $\HH$ is $\Z_2$-graded by an operator
 \be\label{jz} (-1)^F=\begin{pmatrix} I & 0 \cr 0 & -I \end{pmatrix}, \ee
 where $I$ is the identity operator.
 The hamiltonian $H$ commutes with the $\Z_2$-grading, but it is supposed to be the square of an odd self-adjoint  operator $Q$,
 that is, a self-adjoint  operator that anticommutes with $(-1)^F$:
 \be\label{zz} Q=\begin{pmatrix}0 & P\cr P^\dagger & 0\end{pmatrix},~~~~H=Q^2=\begin{pmatrix} P P^\dagger & 0\cr 0 & P^\dagger P\end{pmatrix}.\ee
We then consider a random ensemble for $Q$ defined by the measure $\exp(-N\Tr\,T(Q^2))$, for some suitable function $T$.   If the function $T$
is linear, this is a Gaussian ensemble; here we are interested in a non-Gaussian ensemble.   Ensembles of this type have been introduced (originally in the Gaussian case)
in \cite{W,V,AZ}.

Random supersymmetric quantum mechanics leads easily to the $E^{-1/2}$ behavior of the density of eigenvalues near $E=0$ that we see in eqn. (\ref{nz}).
Let $\mu$ be an eigenvalue of $Q$, and let $f(\mu)\d\mu$ be the corresponding normalized density of eigenvalues.   Since the ensemble for $Q$ is invariant
under $Q\to -Q$, we have $f(\mu)=f(-\mu)$.  The condition  $f(0)\not=0$ is generic.    Now since $H=Q^2$, the relation between an eigenvalue $E$ of $H$ and an eigenvalue $\mu$
of $Q$ is $E=\mu^2$.    Since $f(\mu)\d \mu=f(\sqrt E)\d E/2\sqrt E$, the density of eigenvalues of $H$ behaves as $E^{-1/2}$ near $E=0$.   So an ensemble of this
kind is a good candidate for interpreting the formula (\ref{nz}) for $Z_{\h U}(\beta)$.  

For such a  ``supersymmetric'' matrix ensemble, there is again a version of topological recursion.  Applying this to the spectral curve of eqn. (\ref{wz}), we get
an expansion of $\left\langle \Tr_\HH \,\exp(-\beta H)\right\rangle$   in powers of $e^{-2S}$.  More generally, we can study the expectation value of a product
of traces $\left\langle\prod_{i=1}^n \Tr_\HH\,\exp(-\beta_i H)\right\rangle$.
The terms in the expansion are related to the  volumes of supermoduli spaces in the same manner as in the bosonic case.

In this way, Stanford and I deduced a recursion relation that determines the volumes of the supermoduli spaces $\MM_{\sg,\vec b}$. Moreover, we were 
able to prove this formula by finding analogs of results of Mirzakhani and of Eynard and Orantin.
By imitating Mirzakhani's derivation, we obtained a Mirzakhani-style recursion relation for the volumes of supermoduli spaces.   (See also \cite{Penn} for a related super
McShane identity.)
And similarly to the arguments of Eynard and Orantin, we showed that the recursion relation that comes from the matrix ensemble agrees with the Mirzakhani-style recursion relation.
Thus the relation of random matrices, volumes, and gravity has a perfect counterpart in the supersymmetric case.

Actually, the supersymmetric random matrix ensemble that was just described potentially involves
 an integer invariant, the index of the operator $P$, or equivalently the difference in dimension
between the even and odd subspaces of $\HH$.
    It turns out that to compute volumes of supermoduli spaces  $\MM_{\sg,\vec b}$ of super Riemann surfaces with geodesic boundaries, one should take the index to vanish.   
    But what happens if the index is nonzero?
At least in low orders of the expansion in powers of $e^{-S}$,  the answer is known:  the same type of ensemble based on the same spectral curve but
 with a nonzero index can be used to compute volumes of moduli spaces of super Riemann surfaces with
Ramond punctures as well as geodesic boundaries.   (A Ramond puncture represents a certain type of singularity in the superconformal structure of a super Riemann surface;
roughly, the spin structure is branched over a Ramond puncture.) 
For this, one takes the size of the matrices to infinity, keeping the index fixed. 
A complete proof to all orders in $e^{-S}$ is not yet available in this case.  

Ramond punctures are the only punctures that add something essentially new to the computation of volumes, once one has already analyzed surfaces with geodesic
boundaries.
For purposes of computing volumes, a puncture of an ordinary Riemann surface, or a Neveu-Schwarz puncture of a super Riemann surface
(which is a similar notion),   is equivalent to the small $b$ limit of a geodesic boundary of length $b$.    Nothing like that is true for Ramond punctures.
 
Another important detail concerns the spin structures on a Riemann surface, which may be ``even'' or ``odd.''  The two cases are distinguished by a
topological invariant,  interpreted by Atiyah as the mod 2 index of the Dirac operator \cite{AtiyahModtwo}.  Accordingly, $\MM_\sg$  has two components, parametrizing
super Riemann surfaces with even or odd spin structure.  In the case of a super Riemann surface with boundary,  there is another topological invariant (also
interpretable as a mod 2 index) that specifies
 whether the spin structure on a given boundary component is even or odd.   To get a full answer for the volumes in all cases, one
has to consider also a somewhat different matrix ensemble in which it is still true that $H=Q^2$ (leading again to $\h\rho(E)\sim E^{-1/2}$ near $E=0$),
but a $\Z_2$-grading by $(-1)^F$ is not assumed.

Finally, one can make a similar analyis for unorientable two-manifolds.  
Even if $\Sigma$ is unorientable, we can still define a moduli space of hyperbolic structures on $\Sigma$.
These moduli spaces  are not symplectic,  but JT gravity or supergravity determines  volume forms on them.
It turns out that these volume forms
can be computed using the Reidemeister or Ray-Singer torsion of a flat connection (whose structure group is the full symmetry group of $\H$ or $\h\H$, including
symmetries that reverse the orientation).  
In the bosonic case, the same volume form was first defined by Norbury \cite{N} in another way.
The moduli space volumes are in most cases divergent if $\Sigma$ is unorientable, but the  volume forms can  be compared to random matrix theory.
In a detailed analysis \cite{SW}, all ten standard types of random matrix ensemble \cite{AZ} make an appearance.

\noindent{\it Acknowledgment}  Research  supported in part by  NSF Grant PHY-1911298.

\bibliographystyle{unsrt}

\begin{thebibliography}{99}

\bibitem{AB}
M. F. Atiyah and R. Bott, ``The Yang-Mills Equations Over Riemann Surfaces,'' Phil. Trans. Roy. Soc. Longon {\bf A308} (1982) 523.

\bibitem{Norbury}
P. Norbury, ``A New Cohomology Class On The Moduli Space Of Curves,'' arXiv:1717.03662.

\bibitem{M2}M. Mirzakhani, ``Weil-Petersson Volumes And Intersection Theory On The Moduli Space Of Curves,''
Journal of the American Mathematical Society {\bf 20} (2007) 1-23.

\bibitem{MZ}
Yu. Manin and P. Zograf, ``Invertible Cohomological Field Theories And Weil-Petersson Volumes,'' arXiv:math/9902051.

\bibitem{Z}
P. Zograf, ``On The Large Genus Asymptotics Of Weil-Petersson Volumes,'' arXiv:0812.0544.

\bibitem{M1}
M. Mirzakhani, ``Simple Geodesics  and Weil-Petersson Volumes
Of Moduli Spaces of Bordered Riemann Surfaces,'' Invent. Math. {\bf 167}
(2007) 179-222.

\bibitem{WittenOld}
E. Witten, ``Two-Dimensional Gravity And Intersection Theory On Moduli Space,''
Surveys Diff. Geom. {\bf 1} (1991) 243-310.

\bibitem{BK} E. Brezin and V.  Kazakov, ``Exactly Solvable Field Theories Of Closed Strings,'' Phys. Lett. {\bf B236} (1990) 144.

\bibitem{DS} M.  Douglas and S. Shenker,  ``Strings In Less Than One Dimension,'' Nucl. Phys. {\bf B335}  (1990) 635.

\bibitem{GM} D. J. Gross and A.  Migdal,  ``Nonperturbative Two-Dimensional Quantum Gravity,'' Phys. Rev. Lett. {\bf 64} (1990) 127.

\bibitem{K}M. Kontsevich, ``Intersection Theory On The Moduli Space Of Curves And The Matrix Airy Function,'' Commun. Math. Phys. {\bf 147} (1992)
1-23.



\bibitem{SSS}
P. Saad, S. Shenker, and D. Stanford, ``JT Gravity As A Matrix Integral,'' arXiv:1903.11115

\bibitem{EO}
B. Eynard and N. Orantin, ``Invariants of Algebraic Curves and Topological Expansion,''
Communications in Number Theory and Physics {\bf 1}  (2007) 347-452.




\bibitem{SW}
D. Stanford and E. Witten, ``JT Gravity And The Ensembles Of Random Matrix Theory,''   arXiv:1907.03363.

\bibitem{Ki} A. Kitaev, ``A Simple Model of Quantum Holography.''
\url{http://online.kitp.ucsb.edu/online/entangled15/kitaev/}, \url{http:
//online.kitp.ucsb.edu/online/entangled15/kitaev2/},  talks at KITP, April
7, 2015 and May 27, 2015.

\bibitem{MS}
J. Maldacena and D. Stanford, ``Comments on the Sachdev-Ye-Kitaev Model,'' Phys. Rev. {\bf D94} (2016) 106002.
arXiv:1604.07818.

\bibitem{SW2}
D. Stanford and E. Witten, ``Fermionic Localization of the Schwarzian
Theory,'' arXiv:170-3.04612, JHEP {\bf 10} (2017) 008.

\bibitem{DH}
J. J. Duistermaat and G. J. Heckman, ``On the Variation in the Cohomology of the
Symplectic Form of the Reduced Phase Space'' Inv. Math. {\bf 69} 
(1982) 259Ð268.

\bibitem{AB2}
M. F. Atiyah and R. Bott, ``The Moment Map And Equivariant Cohomology,''
Topology {\bf 23} (1984) 1Ð28.

\bibitem{York}
J. W. York, 
``Role of Conformal Three-Geometry in the Dynamics of Gravitation,''  Phys. Rev. Lett. {\bf 28} (1972)  1082

\bibitem{GH}
G. W. Gibbons and S. W. Hawking, ``Action Integrals and Partition Functions in Quantum Gravity,'' Phys. Rev. {\bf D15} (1977) 2752.

\bibitem{EO2}
 B. Eynard and N. Orantin, ``Invariants of Algebraic Curves and Topological Expansion,''
arXiv:math-ph/0702045.


\bibitem{MaldaEtal}
J. Maldacena, D. Stanford, and Z. Yang,
``Conformal Symmetry And Its Breaking In Two Dimensional Anti-de-Sitter Space,'' arXiv:1606.01857,
Prog. Theor. Exp. Phys.  {\bf 2016}  (2016) 12C104.


\bibitem{Malda}
J. Maldacena, ``The Large $N$ Limit Of Superconformal Field Theories
and Supergravity,'' Int. J. Theor. Phys. {\bf 38} (1999) 1113-1133.

\bibitem{Hawking}
S. W. Hawking, ``The Path Integral Approach To Quantum Gravity,'' in {\it General Relativity - An Einstein Centenary Survey} (Cambridge University Press, 1977), eds. 
S. W. Hawking and W. Israel.

\bibitem{thooft}
G. 't Hooft, ``A Planar Diagram Theory For Strong Interactions,''
Nucl. Phys. {\bf B72} (1974) 461-73.

\bibitem{bpiz}
E. Br\'{e}zin,  C. Itzykson, G. Parisi, and J. B. Zuber, ``Planar Diagrams,'' Comm. Math. Phys. {\bf 59} (1978) 35-51.

\bibitem{DW}
R. Dijkgraaf and E. Witten, ``Developments In Topological Gravity,'' in {\it Topology and Physics}, eds. C. N. Yang, M. L. Ge, and Y. H. He (World Scientific, 2018),
arxiv:1804.03275.

\bibitem{W}
J. Wishart, ``The Generalized Product Distribution In Samples From A Normal Multivariate Population,'' Biometrika {\bf 20A} (1928) 32-52.

\bibitem{V}
J. J. M. Verbaarschot and I. Zahed,
``Spectral Density of the QCD Dirac Operator
Near Zero Virtuality,'' Phys. Rev. Let. {\bf 70} (1993) 3853-3855, arXiv:hep-th/9303012.

\bibitem{AZ}
A. Altland and M. R. Zirnbauer, ``Nonstandard Symmetry Classes in Mesoscopic Normal-Superconducting
Hybrid Structures,'' Phys. Rev. {\bf B55} (1997) 1142.

\bibitem{Penn}
Yi Huang, R. C. Penner, and A. Zeitlin,
``Super McShane Identity,'' arXiv:1907.09978.

\bibitem{AtiyahModtwo}
M. F. Atiyah, ``Riemann Surfaces and Spin Structures,'' Ann. Sci. de l'\'{E}.N.S. {\bf 4} (1971) 47-62.

\bibitem{N}
P. Norbury, ``Lengths of Geodesics On Non-Orientable Hyperbolic Surfaces,'' Geom. Dedicatae {\bf 134} (2008) 153-76, arXiv:math/0612128.


\end{thebibliography}

\end{document}